\documentstyle[12pt]{article}

\newcommand{\sect}[1]{\setcounter{equation}{0}\section{#1}}
\newcommand{\subsect}[1]{\subsection{#1}}

\def\be{\begin{equation}}
\def\ee{\end{equation}}
\def\bea{\begin{eqnarray}}
\def\eea{\end{eqnarray}}

\def\pois#1#2{\left\{ #1,#2 \right\}}
\def\te{\chi}
\def\de{\theta}
\def\h{\hbar}

\def\vvv{v}

\def\>#1{{\bf #1}}
\def\1{\'{\i}}
\def\R{\rm I\kern-.2em R}

\begin{document}

\thispagestyle{empty}

\hfill\today

\vspace{2cm}

\begin{center}
{\LARGE{\bf{Non-standard quantum (1+1)}}}

{\LARGE{\bf{Poincar\'e group:}}}

{\LARGE{\bf{a $T$--matrix approach}}}

\end{center}

\bigskip\bigskip\bigskip

\begin{center}
A. Ballesteros$^{1}$, F.J. Herranz$^{1}$, M.A. del
Olmo$^{2}$, C.M. Pere\~na$^{3}$ and M. Santander$^{2}$
\end{center}

\begin{center}
\em

{ { {}$^{1}$ Departamento de F\1sica, Universidad de Burgos}
\\  E-09003, Burgos, Spain}

{ { {}$^{2}$Departamento de F\1sica Te\'orica, Universidad de Valladolid }\\
E-47011, Valladolid, Spain }

{ $^{3}$ Departamento de F\1sica Te\'orica, Universidad Complutense\\
    E-28040, Madrid, Spain}

\end{center}
\rm

\bigskip

\bigskip\bigskip\bigskip

\begin{abstract}

The Hopf algebra dual form for the non--standard uniparametric deformation of
the (1+1) Poincar\'e algebra $iso(1,1)$ is deduced. In this framework, the
quantum coordinates that generate $Fun_w(ISO(1,1))$ define an infinite
dimensional Lie algebra. A change in the basis of the dual form is obtained
in order to compare this deformation to the standard one. Finally, a
non--standard quantum Heisenberg group acting on a quantum Galilean plane is
obtained.

\end{abstract}

\newpage

\sect{Introduction}

Quantum deformations were introduced  as quantizations of some integrable
models characterized by quadratic Poisson brackets \cite{G}. The
essential relation between these brackets and Lie bialgebras was established
by Drinfel'd \cite{DrN}. This link provides a framework to understand the
quantization of such Poisson--Lie structures as a dual process to the
(bialgebra) deformation of universal enveloping Lie algebras \cite{Tak}. More
recently, the $T$--matrix approach has been introduced in such a way that the
Hopf algebra dual form $T$ summarizes   duality in a universal
(representation independent) setting \cite{FG,FGT,F}. The transfer
matrices of certain integrable systems can be seen as particular realizations
of the dual form $T$ and, conversely, the construction of new integrable
models could be guided by the obtention of new $T$--matrices. In general, this
approach provides a ``canonical" formalism in which quantum objects are
defined in a completely similar language to their classical counterparts.
This ``proper" setting provides a generalization of many group theoretical
results to the non--commutative cases in a straightforward --although
possibly cumbersome-- way (for instance, see \cite{BCiGST} for an example in
the context of $q$--special function theory). From an algebraic point
of view, the $T$--matrix approach emphasizes the equivalent role that both
quantum algebras and quantum groups play as Hopf algebra deformations and
reveals the importance of solvable Lie algebras in this context, a fact that
has been also pointed out in \cite{Lyak}.

The Lie bialgebra structures compatible with a determined Lie algebra
give a primary characterization for its quantum algebra deformations. So far,
the deformations coming from the (non--degenerate) coboundary Lie bialgebra
structures classified by Belavin and Drinfel'd \cite{BD} have been deeply
studied (the so--called ``standard" deformations \cite{Dr,Ji}). In some cases,
their corresponding $T$--matrices have been deduced and, by using contraction
methods, these results have been extended to some quantum non--semisimple
groups \cite{FG,FGT,F,BCGP,CJ}.

In this paper we deal with the non--semisimple inhomogeneous algebra
$iso(1,1)\simeq t_2\odot so(1,1)$ with classical commutation rules
\be
[K ,P_\pm]=\pm 2P_\pm ,\qquad [P_+,P_-]= 0.\label{bc}
\ee
As a real form, the algebra $iso(1,1)$ is isomorphic to the (1+1) dimensional
Poincar\'e algebra: $K$ and $P_\pm$ generate, respectively, the boosts and the
translations along the light-cone. A non--standard coboundary Lie bialgebra
$(iso(1,1),\delta^{(n)})$ is generated by the classical $r$--matrix
\be
r^{(n)}=K \wedge P_+.\label{bd}
\ee
As usual, $\delta^{(n)}(X)=[1\otimes X + X\otimes 1,r^{(n)}]$ and  $r^{(n)}$
verifies the Classical Yang--Baxter Equation (CYBE). The set of coboundary
structures for this algebra is completed by the Lie bialgebra
$(iso(1,1),\delta^{(s)})$ given by
\be
r^{(s)}=K \wedge (P_- + P_+),\label{bdd}
\ee
that generates the standard deformation \cite{BCGpho,PL} (this $r^{(s)}$
fulfills the modified CYBE). As a particular feature of the (1+1) Poincar\'e
algebra \cite{LBC,Zakr}, there also exists a non--coboundary $iso(1,1)$
bialgebra with cocommutator
\be
\delta^{(nc)}(K)=0,\qquad \delta^{(nc)}(P_\pm)=P_\pm\wedge K.
\label{bde}
\ee
It can be easily checked that the quantum Poincar\'e algebra of \cite{VK} is a
deformation of $iso(1,1)$ ``in the direction" of $\delta^{(nc)}$.

Starting from the quantization $U_w iso(1,1)$ (reviewed in Section 2) of the
Lie bialgebra $(iso(1,1),\delta^{(n)})$, the construction of its Hopf algebra
dual form is developed in Section 3. The main result of this
section is that the non--standard quantum $T$--matrix is written as a product
of usual exponentials.  We recall that the (1+1) Poincar\'e dual form given
in \cite{BCGP} (that corresponds to the quantum algebra \cite{VK}) does
contain $q$--exponential factors.

In Section 4, the non--standard quantum $ISO(1,1)$ group
is deduced. We emphasize a relevant difference of
the result so obtained with respect to the $T$--matrices already known: in
\cite{FG,FGT,BCGP,CJ} the classical dynamical variables (i.e., the group
coordinates under a certain factorization of the group elements) generate a
(solvable) finite dimensional Lie bialgebra $(g_x,\delta)$, and the quantum
group is constructed as a deformation $U_q(g_x)$. In our case, the Poisson
algebra $g_x$ defined by the Sklyanin bracket coming from (\ref{bd}) is an
infinite dimensional Lie algebra, and its quantization is obtained by
applying the Weyl prescription onto this infinite--dimensional object. In
particular, this enhances the differences between the $T$--matrix of
\cite{BCGP} and the one presented here.

On the other hand, the symmetric character of the
$T$--matrix is used to transform in a consistent way the light--cone quantum
coordinates into the quantum  space--time cartesian coordinates. This allows
us to compare the results here obtained to those of \cite{PL} that correspond
to a ``standard" quantum $ISO(1,1)$ group dual to the deformation generated by
(\ref{bdd}). By following the scheme developed in \cite{PL}, it is possible to
understand $U_w iso(1,1)$ as a quantization of the Poisson--Lie similarity
group $SB(2)$.

To end with, Section 5 introduces a new quantum Heisenberg group obtained
from the previous results by means of a formal transformation involving
Clifford dual units. This quantum group acts naturally on a quantum Galilei
plane, and can be easily embedded within the classification given in
\cite{HLR}.


\sect{Non--standard quantum $iso(1,1)$ algebra}

It is well known that in $sl(2,\R)$ there exist two non trivial  (coboundary)
Lie bialgebra structures \cite{PC}. The first one is generated by
$r^{(s)}=J_+\wedge J_-$ and underlies the (standard) Drinfel'd--Jimbo
deformation for this algebra. The quantization of the non--standard bialgebra
(with $r^{(n)}=J_3\wedge J_+$) was developed in \cite{Ohn}.

At a classical level, the involutive automorphism of
$sl(2,\R)$ given by $S(J_3,J_\pm)=(J_3,-J_\pm)$ induces an In\"on\"u--Wigner
contraction of this algebra that is obtained as the limit $\varepsilon\to 0$
of the transformation $(K,P_\pm):=\Gamma(J_3,J_\pm)=(J_3,\varepsilon J_\pm)$.
The algebra that arises under such a contraction is just (\ref{bc}). This
procedure can be extended to the quantum case by taking into account the
following ``quantum" automorphism of $U_z sl(2,\R)$
\be
S_q(J_3,J_\pm;z)=(J_3,-J_\pm;-z).\label{cc}
\ee
Involution $S_q$
gives rise to a generalized In\"on\"u--Wigner contraction:
\be
(K,P_\pm;w):=\Gamma_q(J_3,J_\pm;z)=(J_3,\varepsilon
J_\pm;z/\varepsilon),\label{cd}
\ee
where $K$ and $P_\pm$ and
$w$ are the Lie generators and the deformation parameter of the contracted
quantum algebra $U_w iso(1,1)$. By applying the transformation $\Gamma_q$ onto
the deformation given in \cite{Ohn} and making the limit $\varepsilon\to 0$,
the quantum algebra $U_w iso(1,1)$ is obtained (see \cite{Bey}). It is worth
recalling that this quantum algebra was firstly discovered in \cite{LM} with
no reference to contraction procedures.

\noindent {\bf Proposition 2.1.}
\begin{em}
Let $K$ and $P_\pm$ the infinitesimal generators of the Lie algebra
$iso(1,1)$. The Hopf algebra $U_w iso(1,1)$ is given by the coproduct,
counit, antipode
\bea
&&\Delta P_+ =1 \otimes P_+ + P_+\otimes 1,\nonumber\\
&&\Delta P_- =e^{-wP_+} \otimes P_- + P_-\otimes e^{wP_+},\label{ce}\\
&&\Delta K =e^{-wP_+} \otimes K + K\otimes e^{wP_+},\nonumber
\eea
\be
\epsilon(X) =0;\qquad
\gamma(X)=-e^{wP_+}\ X\ e^{-wP_+},\ \ \mbox{for $X\in \{K,P_\pm\}$},
\label{cf}
\ee
\be
\gamma(P_+)=-P_+,\quad  \gamma(P_-)=-P_- ,\quad
\gamma(K)=-K+2\sinh (wP_+);\label{cg}
\ee
and the commutation relations
\be
[K,P_+]=2\,\frac{\sinh (wP_+)}w,\quad [K,P_-]=- 2\, P_-\cosh (wP_+) ,
\quad  [P_+,P_-]= 0.\label{ch}
\ee
This structure is a quantization of the non--standard Lie
bialgebra of $iso(1,1)$ generated by $r^{(n)}= K \wedge P_+$.
\end{em}

The center of  $U_w iso(1,1)$ is generated by
\be
C_w= 2 P_-\frac{\sinh (wP_+)}w .\label{ci}
\ee
A  classical $3\times 3$ matrix representation of $iso(1,1)$ is given by:
\be
D(K)=\left(\begin{array}{ccc}
 0 & 0 & 0 \\
 0 & 0 & - 2\\
 0 & -2 & 0  \end{array}\right),\
D(P_+)=\left(\begin{array}{ccc}
 0 & 0 & 0 \\
 1 & 0 & 0\\
 -1 & 0 & 0  \end{array}\right),\
D(P_-)=\left(\begin{array}{ccc}
 0 & 0 & 0 \\
 1 & 0 & 0\\
 1 & 0 & 0  \end{array}\right),\label{cj}
\ee
and a quantum matrix realization of $U_w iso(1,1)$ coincides with the
classical one:
\be
D_q(X)=D(X),\quad X\in\{K,P_\pm\}.\label{ck}
\ee


\sect{The universal $T$--matrix}

Following \cite{FG,FGT} we consider the Hopf algebra dual form
\be
T=\sum_{abc} p_{abc}\otimes X^{abc} ,
\label{dc}
\ee
where $X^{abc}$ is a basis for $U_w iso(1,1)$ and its dual basis
 $p_{abc}$ generates $Fun_w (ISO(1,1))$
(sumation over repeated indices will be omitted from now on).

The $T$--matrix is constructed here starting from $U_w iso(1,1)$. A suitable
basis to work   is $X^{abc}= A_-^a {A_+}^b H^c$ (hereafter this ordering
will be preserved), where
\be
{A_+}=P_+,\qquad {A_-}=e^{-wP_+}P_-,\qquad H=e^{wP_+}K.
\label{dd}
\ee
The following coproduct and commutation rules are deduced from (\ref{ce})
and (\ref{ch}):
\bea
&&\Delta ({A_+}) =1 \otimes {A_+} + {A_+}\otimes
1, \nonumber\\
&&\Delta ({A_-}) =e^{-2w{A_+}} \otimes {A_-} + {A_-}\otimes 1,
\label{df}\\
&&\Delta (H) =1 \otimes H + H\otimes e^{2w{A_+}}
\nonumber; \eea
\be
[H,{A_+}]=\frac{e^{2w{A_+}}-1}w, \quad [H,{A_-}]=-
2\, {A_-}e^{2w{A_+}} , \quad  [{A_+},{A_-}]=0 \label{dg}.
\ee

The dual basis to $X^{abc}$ is defined by
\be
\langle p_{abc}, X^{lmn} \rangle =\delta_a^l\,\delta_b^m\,\delta_c^n.
\label{dh}
\ee
and we set ${\hat a_-}=p_{100},{\hat a_+}=p_{010}$ and ${\hat \te}=p_{001}$.

\noindent {\bf Theorem 3.1.}
\begin{em} The dual basis $p_{qrs}$ can be expressed in terms of the dual
coordinates ${\hat a_-},{\hat a_+},{\hat \te}$ in the form
\be
p_{qrs}=\frac{{\hat a_-}^q}{q!}\frac{{\hat a_+}^r}{r!}\frac{{\hat
\te}^s}{s!}, \label{di}
\ee
and the Hopf algebra dual form $T$ reads
\be
T=e^{{\hat a_-} {A_-}}e^{{\hat a_+} {A_+}}e^{{\hat \te} H}.
\label{dj}
\ee
\end{em}

The proof of this statement follows by considering the structure tensor
$F$ that gives us the coproduct of an arbitrary element of $U_w
iso(1,1)$
\be
\Delta(X^{abc}):=F_{lmn;qrs}^{abc}\,X^{lmn}\otimes X^{qrs}
\label{dk}
\ee
and, by means of duality, the product in $Fun_w (ISO(1,1))$
\be
p_{lmn}\, p_{qrs}=F_{lmn;qrs}^{abc}\,p_{abc}.
\label{dl}
\ee
In particular, it is easy to see that
\bea
&& F_{000;qrs}^{abc}=\delta_q^a\,\delta_r^b\,\delta_s^c,\nonumber\\
&& F_{lmn;000}^{abc}=\delta_l^a\,\delta_m^b\,\delta_n^c,\label{dm}\\
&&
F_{lmn;qrs}^{000}=\delta_l^0\,\delta_m^0\,\delta_n^0
\,\delta_q^0\,\delta_r^0\,\delta_s^0.\nonumber
\eea

Two recurrence relations for $F$ can be deduced by observing that
\bea
&&\Delta(X^{abc})=\Delta ({A_-})\Delta(X^{(a-1)bc}),\label{dn}\\
&&\Delta(X^{abc})=\Delta ({A_+})\Delta(X^{a(b-1)c}),\label{do}
\eea
and using (\ref{df}) and (\ref{dg}). In this way we find:
\bea
&& F_{lmn;qrs}^{abc}=\sum_{k=0}^m F_{lkn;(q-1)rs}^{(a-1)bc}
\frac{(-2w)^{m-k}}{(m-k)!} + F_{(l-1)mn;qrs}^{(a-1)bc},\quad
a,q,l\geq1,\label{dp}\\ && F_{lmn;qrs}^{abc}=F_{l(m-1)n;qrs}^{a(b-1)c}
+ F_{lmn;q(r-1)s}^{a(b-1)c},\quad b,m,r\geq1.
\label{dq}
\eea
The recurrence relation corresponding to $\Delta(X^{abc})=
\Delta(X^{abc-1})\Delta (H)$ is much harder to find in general. However,
for our purposes, we shall only need some particular cases of it.

The relations (\ref{dp}) y (\ref{dq}), together with (\ref{dm}), lead to
\bea
&&
F_{100;qrs}^{abc}=a\,\delta_{q+1}^a\,\delta_{r}^b\,\delta_s^c,\qquad
a\geq 1,\label{dr}\\ &&
F_{010;qrs}^{abc}=b\,\delta_q^a\,\delta_{r+1}^b\,\delta_s^c,\qquad
b\geq 1,\label{ds}
\eea
and by simple considerations we can also find  that
\be
F_{lmn;001}^{abc}=c\,\delta_l^a\,\delta_m^b\,\delta_{n+1}^c,\quad
c\geq 1.
\label{dt}
\ee
These are the elements of $F$ that are relevant in order to compute the dual
basis. Now, with the aid of  (\ref{dk}) and (\ref{dr}), the following relation
holds \be
p_{100}\, p_{(q-1)rs}=F_{100;(q-1)rs}^{abc}\,p_{abc}=a\,
\delta_{q}^a\,\delta_{r}^b\,\delta_s^c \,p_{abc} =q\,p_{qrs},
\label{du}
\ee
so
\be
p_{qrs}=\frac{{\hat a_-}}{q}\,p_{(q-1)rs}=
\dots=\frac{{\hat a_-}^q}{q!}\,p_{0rs}.
\label{dv}
\ee
Straightforward calculations based on (\ref{ds}) and (\ref{dt})
together with the fact that $p_{000}=1$ complete the proof of the theorem.


\sect{The quantum group $Fun_w (ISO(1,1))$}

The structure tensor $F$ also allows us to deduce the commutation rules
between the generators of $Fun_w (ISO(1,1))$. As a result, we enunciate the
following proposition:

\noindent {\bf Proposition 4.1.}
\begin{em} The dual coordinates ${\hat a_-}$, ${\hat a_+}$ and ${\hat \te}$
satisfy \be
[{{\hat \te}},{{\hat a_+}}]=w(e^{ 2{\hat \te}} -1) ,\quad
[{{\hat \te}},{{\hat a_-}}]=0,\quad [{{\hat a_+}},{{\hat a_-}}]=- 2w  {\hat
a_-} . \label{ec}
\ee
\end{em}

\noindent \begin{em} Proof: \end{em} The commutation relation of any two
elements of $Fun_w (ISO(1,1))$ is
\be
[p_{lmn},p_{qrs}]=(F_{lmn;qrs}^{abc}-F_{qrs;lmn}^{abc})p_{abc} .
\label{ed}
\ee
The explicit expressions (\ref{ec}) for  $[{{\hat \te}},{{\hat
a_-}}]\equiv[p_{001},p_{100}]$ and $[{{\hat a_+}},{{\hat
a_-}}]\equiv[p_{010},p_{100}]$ are straightforwardly derived from the
relations involving $F$ in the previous section. In spite of the absence of a
third general recurrence relation, we can again find the particular values of
$F$ involved in  \be
[{{\hat \te}},{{\hat a_+}}]=(F_{001;010}^{abc}-F_{010;001}^{abc})p_{abc}.
\label{ef}
\ee
{}From (\ref{ds}) we obtain that $F_{010;001}^{abc}=
b\,\delta_0^a\,\delta_{1}^b\,\delta_1^c$. On the other hand,
$F_{001;010}^{abc}$ gives the coefficient of the term $H\otimes A_+$
 in $\Delta(X^{abc})$. Such a term appears either when $b=c=1$ (this
contribution annihilates the previous $p_{011}$ term) or in the cases
$a=b=0$ and $c$ arbitrary (note that, since $[H,A_\pm]$ does not produce
$H$, this generator cannot appear as a byproduct of reordering processes). As
a consequence,
\be
F_{001;010}^{abc}\, p_{abc}=p_{011}+ 2w\,p_{001}
+4w\,p_{002}+\dots+2^k w\,p_{00k}+\dots
\label{efd}
\ee
Therefore,
\be
[{{\hat \te}},{{\hat a_+}}]=w\sum_{k=1}^\infty
\frac{1}{k!} 2^k{\hat \te}^k=w(e^{ 2{\hat \te}} -1).
\label{efe}
\ee
Note that two algebras (\ref{ec}) with parameters $w$ and $w'$ (both
of them different from zero) are isomorphic.

Since the classical fundamental representation of $iso(1,1)$
(\ref{cj}) is a fundamental one for the quantum algebra $U_w iso(1,1)$, the
specialization of the $T$--matrix (\ref{dj}) to this realization is just the
classical $ISO(1,1)$ group element with non--commutative entries:
\be
T^{D_q}=
e^{{\hat a_-} D_q ({A_-}) }e^{{\hat a_+}   D_q ({A_+}) }e^{{\hat \te}
D_q (H) }=\left(\begin{array}{ccc}
 1 & 0 & 0 \\
 {\hat a_-} + {\hat a_+} & \cosh 2{\hat \te} & -  \sinh 2{\hat \te}\\
 {\hat a_-} - {\hat a_+} & -\sinh 2{\hat \te} & \cosh 2{\hat \te}
\end{array}\right).
\label{eg}
\ee
Thus, the multiplicative property of $T$ provides the coproduct for $Fun_w
(ISO(1,1))$ (note that $D_q ({P_\pm})=D_q ({A_\pm}),D_q ({K})=D_q (H)$).

\noindent {\bf Theorem 4.2.}
\begin{em} The Hopf algebra  $Fun_w (ISO(1,1))$ is given by the commutation
rules (\ref{ec}), coproduct
\bea
&&\Delta({\hat \te})={\hat \te}\otimes 1 + 1 \otimes {\hat \te},\cr
&&\Delta({\hat a_+})={\hat a_+}\otimes 1 + \cosh 2{\hat \te}\otimes  {\hat
a_+} + \sinh 2{\hat \te}\otimes {\hat a_+} ,\label{eh}\\
&&\Delta({\hat a_-})={\hat a_-}\otimes 1 + \cosh 2{\hat \te}\otimes  {\hat
a_-} - \sinh 2{\hat \te}\otimes {\hat a_-} ;\nonumber
\eea
counit and antipode
\be
\epsilon(X) =0,\qquad
X\in\{   {\hat a_+},   {\hat a_-}, {\hat \te}\};
\label{ei}
\ee
\be
\gamma(\hat  \chi)=-\hat  \chi,\quad
 \gamma( {\hat a_+})=- e^{-2{\hat \te}}{\hat a_+},\quad
\gamma( {\hat a_-})=-e^{2{\hat \te}}{\hat a_-} .
\label{ej}
\ee
\end{em}


The ``coincidence" between classical and quantum coordinates can be also
extracted from the Poisson--Lie structure underlying this quantization.
A Poisson--Hopf algebra of smooth functions on the group $ISO(1,1)$ can be
constructed by means of the classical non--standard $r$--matrix
$r^{(n)}=K\wedge P_+$. The Sklyanin bracket induced from this $r$--matrix is
\bea
&&\pois{f}{g}=r^{\alpha\beta}\left(X_{\alpha}^L f X_{\beta}^L g -
X_{\alpha}^R f
X_{\beta}^R g \right) \cr
&&=   X_{K}^L f X_{P_+}^L g -  X_{P_+}^L f X_{K}^L g  -
X_{K}^R f X_{P_+}^R g
+ X_{P_+}^R f X_{K}^R g . \label{ek}
\eea
We use the
classical matrix representation similar to (\ref{eg}) to obtain the left
and right invariant vector fields:
\be
X_{K}^L=\partial_\te,\qquad
X_{P_+}^L=e^{2\te}\partial_{a_+},\qquad
X_{P_-}^L=e^{-2\te}\partial_{a_-},\label{el}
\ee
\be
X_{K}^R=\partial_\te + 2 a_+ \partial_{a_+}- 2 a_-\partial_{a_-},\qquad
X_{P_+}^R=\partial_{a_+},\qquad
X_{P_-}^R=\partial_{a_-},\label{en}
\ee
and we get the following result:

\noindent {\bf Proposition 4.3.}
\begin{em} The Poisson bracket
\be
\pois{f}{g}=m\circ ( (e^{2\te} -1)
\partial_\te\wedge\partial_{a_+} -
2a_- \partial_{a_+}\wedge\partial_{a_-} )(f\otimes g),
\label{eeo}
\ee
gives the structure of
a Poisson--Hopf algebra to $Fun(ISO(1,1))$ ($m(a\otimes b)=ab$).
\end{em}


In particular for the coordinates $a_+$, $a_-$ and $\te$ we get
\be
\{ {{\hat \te}},{{\hat a_+}}\}=w(e^{ 2{\hat \te}} -1) ,\quad
\{{{\hat \te}},{{\hat a_-}}\}=0,\quad
 \{{{\hat a_+}},{{\hat a_-}}\}=- 2w {\hat a_-} .
\label{eeec}
\ee
This means that the commutation rules (\ref{ec}) can be seen as a Weyl
quantization $\pois{\,}{\,}\rightarrow w^{-1}[\, ,\,]$ of the fundamental
Poisson brackets given by (\ref{eeec}) and that the classical coproduct  does
not change under quantization (compare to \cite{FG,BCGP}).

On the other hand, the FRT \cite{FRT} prescription is consistent with
Theorem 4.2. The (quadratic) commutation rules between the entries of the
quantum matrix (\ref{eg}) that are derived from (\ref{ec}) coincide with the
relations obtained via the FRT procedure and with the aid of the quantum
$R$--matrix  \be
R=1\otimes 1 + w D_q(H)\wedge D_q({A_+}).
\label{ep}
\ee

 Finally, it is worth recalling that in  \cite{FG,FGT,BCGP,CJ}
the $T$--matrix coordinates that generate $Fun_q(G)$ close a solvable
finite dimensional Lie (super) algebra (with the exception of the ``esoteric"
quantum $GL(n)$ \cite{FGT}). In particular, for the $iso(1,1)$ case studied in
\cite{BCGP} it is shown that \be
[\pi,\pi_+]=-z\,\pi_+ ,\quad
[\pi,\pi_-]=-z\,\pi_- ,\quad [\pi_+,\pi_-]=0.
\label{epd}
\ee
This is no longer the case for the
non--standard deformation (\ref{ec}). In fact, if we consider ${\hat
a_+},{\hat a_-}$ and $\hat f_m:=e^{2m{\hat \te}}\, (m\in Z)$, we can identify
$Fun_w(ISO(1,1))$ with an infinite dimensional Lie algebra endowed with
a Hopf algebra structure:
\bea
&&\Delta({\hat f_m})={\hat f_m}\otimes \hat f_m ,\cr
&&\Delta({\hat a_+})={\hat a_+}\otimes 1 + \hat f_1\otimes {\hat
a_+},\label{epe}\\
&&\Delta({\hat a_-})={\hat a_-}\otimes 1 + \hat f_{-1}\otimes {\hat
a_-} .\nonumber
\eea
Moreover, since
\be
[\hat f_m,{\hat a_+}]=2\,w\,m\,(\hat f_{m+1}-\hat f_m),
\label{eq}
\ee
we have an infinite--dimensional Hopf subalgebra generated by
$\hat f_m$ and $\hat a_+$ (note that $H$ and $A_+$ do generate a Hopf
subalgebra aswell). It would be interesting to know whether this is a general
feature of the non--standard deformations and to investigate the dependence
of the dimensionality of the quantum algebra of coordinates with the
parametrization. In principle, this infinite--dimensional aspect makes the
approach formally closer to the realistic integrable models obtained without
truncating the spectral parameter \cite{FGT}.

\subsect {Standard versus non--standard deformations}

In order to compare the results so far obtained to the ones given in
\cite{PL} it is necessary to perform a ``change of quantum coordinates"
relating the ``light cone" ones $({\hat a_+},{\hat a_-})$ with the time and
space quantum translations $(\hat a_1,\hat a_2)$. We recall the defining
relations for the (1+1) quantum Poincar\'e group described in \cite{PL}. The
group element was
\be
G=\left(\begin{array}{ccc}
 1 & 0 & 0 \\ \hat a_1 & \cosh\hat\de & \sinh\hat\de\\
\hat a_2 &\sinh\hat\de & \cosh\hat\de  \end{array}\right),
\label{er}
\ee
and the commutation rules between the coordinates read
\bea
&& [\hat\de,\hat a_1]=w'\,(\cosh \hat \de -1),\nonumber\\
&& [\hat\de,\hat a_2]= w'\,\sinh \hat \de,\label{es} \\
&& [\hat a_1,\hat a_2]= w'\, \hat a_1. \nonumber
\eea
This  algebra  can
be also written in terms  of   generators  $\hat g_n=e^{n \de}$.
 The standard quantum (1+1) Poincar\'e plane of coordinates  $(\hat
x_1^{(s)},\hat x_2^{(s)})$ characterized by  $[\hat x_1^{(s)},\hat
x_2^{(s)}] = w'\, x_1^{(s)}$ is easily derived from these expressions.

Both   standard and non--standard quantizations turn out to be Weyl
quantizations of the classical coordinates preserving the multiplicative
property of  $T^{D_q}$ and $G$  (cfr. (\ref{eeo}) and eq.
(3.1) in \cite{PL}). This fact suggests the use of the classical relation
between coordinates  as an Ansatz for the change of quantum coordinates:
\be
\hat \de=- 2 {\hat \te},\quad \hat a_1={\hat a_-} + {\hat a_+},
\quad \hat a_2= {\hat a_-}  - {\hat a_+}.
\label{et}
\ee
This redefinition implies that $T^{D_q}$ and $G$ become identical. Hence, we
have two different sets of commutation rules compatible with the
same coproduct $\Delta(G)=G\dot\otimes G$. In particular, the change
(\ref{et}) on (\ref{ec}) provides the non--standard brackets $(w'=-2 w)$
\bea
&& [\hat\de,\hat a_1]=w'\,(\cosh \hat \de -1 - \sinh \de),\nonumber\\
&& [\hat\de,\hat a_2]=w'\,( \sinh \hat \de - (\cosh \de -
1)),\label{eu} \\
&& [\hat a_1,\hat a_2]=w'\,(\hat a_1+ \hat a_2). \nonumber
\eea
The non--standard quantum Poincar\'e plane (whose relations are invariant
under the coaction defined by $G$ on a vector $(x_1^{(n)},x_2^{(n)})$) is
characterized by
\be
[\hat x_1^{(n)},\hat x_2^{(n)}] = w'\, ( x_1^{(n)} + x_2^{(n)}).
\label{ev}
\ee

A comparison between both quantizations shows that, in this quantum basis,
the non--standard deformation seems to be constructed by adding some
additional terms on the standard relations. It is also worth remarking the
absolutely symmetrical role that both coordinates play in the
non--standard case (see the quantum plane relation (\ref{ev})). This kind of
``symmetrical" quantization has been already related to the non--standard
(2+1) deformations at a quantum algebra level \cite{Bey}. In the last Section
we shall analyze this symmetry from a more geometrical point of view.

\subsect {A realization of $T$ with entries in $U_w(iso(1,1))$ }

In general, given a Lie group $G$ it is rather difficult to control  how a
given change of basis in $Fun_q(G)$ can be paired to a transformation of the
generators of $U_q(g)$, $g$=Lie($G$), in such a way that duality is preserved.
In the sequel we see how the solution to this problem is naturally contained
within the $T$--matrix approach.

Let us consider the change of basis (\ref{et}). We are interested in finding
a set $\{A_1,A_2,A_{12}\}$ of generators of $U_w(iso(1,1))$ such that the Hopf
algebra dual form $T$ is preserved:
\be
T=e^{{\hat a_-} {A_-}}e^{{\hat a_+} {A_+}}e^{{\hat \te} H}=
e^{{\hat a_1} {A_1}}e^{{\hat a_2} {A_2}}e^{{\hat \de} A_{12}}.
\label{ex}
\ee

Since ``universal" computations to relate both $U_w(iso(1,1))$ bases are
extremely cumbersome, we can specialize the $T$ matrix to a (fundamental)
representation $Q$ of $Fun_w(ISO(1,1))$ given by
\be
Q(\hat \te)=\left(\begin{array}{cccc}
 0 & 0 & 1 & 0 \\
 0 & 0 & 0 & 0 \\
 0 & 0 & 0 & 0 \\
 0 & 0 & 0 & 0   \end{array}\right),\
Q(\hat a_-)=\left(\begin{array}{cccc}
 0 & 0 & 0 & 0 \\
 0 & 0 & 0 & 1 \\
 0 & 0 & 0 & 0 \\
 0 & 0 & 0 & 0   \end{array}\right),\
Q(\hat a_+)=\left(\begin{array}{cccc}
 0 & 0 & 0 & 1 \\
 0 & -2w & 0 & 0 \\
 0 & 0 & 2w & 0 \\
 0 & 0 & 0 & 0   \end{array}\right).\
\label{ey}
\ee
A straightforward computation shows that
\be
T^Q:= e^{Q(\hat a_-) {A_-}}e^{Q(\hat a_+) {A_+}}e^{Q(\hat \te) H}=
\left(\begin{array}{cccc}
 1 & 0 & H & A_+ \\
 0 & e^{-2w {A_+}} & 0 & A_- \\
 0 & 0 & e^{2w {A_+}} & 0 \\
 0 & 0 & 0 & 1   \end{array}\right).\
\label{ez}
\ee
{}From this point of view, $U_w(iso(1,1))$ is a quantum
similarity group (it is easy to check that $Q(\hat a_+),Q(\hat a_-)$ and
$Q(\hat \te)$ close a similarity algebra $sb(2)$) with non commuting
coordinates $A_+,A_-$ and $H$. Moreover, the coproduct (\ref{df}) is
reproduced by the multiplicative property $\Delta(T^Q)=T^Q\dot\otimes T^Q$.

The new quantum coordinates admit a representation $Q(\hat a_1)$, $Q(\hat
a_2)$ and $Q(\hat \de)$ derived from (\ref{et}) and (\ref{ey}).  By computing
the corresponding exponentials we obtain a second expression for the dual form
$T^Q:= e^{Q(\hat a_1) {A_1}}e^{Q(\hat a_2) {A_2}}e^{Q(\hat \de) A_{12}}$:
\be
T^Q=
\left(\begin{array}{cccc}
 1 & 0 & -2A_{12} & A_1 - A_2 \\
 0 & e^{-2w (A_1 - A_2)} & 0 & \frac{1}{2w}(e^{-2w (A_1 - A_2)} - 2e^{-2w
A_1} + 1)\\
 0 & 0 & e^{2w (A_1 - A_2)} & 0 \\
 0 & 0 & 0 & 1   \end{array}\right).\
\label{ezc}
\ee
The relation between the two $U_w(iso(1,1))$ bases is now clear
\bea
&& A_+= A_1 - A_2,\nonumber\\
&& A_-= \frac{1}{2w}(e^{-2w (A_1 - A_2)} - 2e^{-2w A_1} +
1),\label{ezd}\\
&& H=-2A_{12}.\nonumber
\eea
Note that $\lim_{w\to 0} A_-= A_1 + A_2$ and we recover the usual
classical change of basis.

\subsect {  A Poisson--Lie realization of $U_w(iso(1,1))$ }

The quantum algebra $U_w(iso(1,1))$ can be considered as a Lie bialgebra
deformation of $iso(1,1)$. The cocommutator $\delta^{(n)}$ is a  coboundary
that gives the first order term in $w$ of the antisymmetrized part of the
coproduct (\ref{df}):
\bea
&& \delta^{(n)}(A_+)=0,\nonumber\\
&& \delta^{(n)}(A_-)=2w \,A_-\wedge A_+,\label{eze}\\
&& \delta^{(n)}(H)=2w \,H\wedge A_+.\nonumber
\eea

As we mentioned above, if we take into account the representation  (\ref{ey})
of the Hopf algebra $Fun_w(ISO(1,1))$, the commutation rules (\ref{ec}) are
linearized into the Lie algebra $sb(2)$
\be
[Q({{\hat \te}}),Q({{\hat a_+}})]=2w\,Q({\hat \te}) ,\quad
[{{\hat \te}},{{\hat a_-}}]=0,\quad [Q({{\hat a_+}}),Q({{\hat a_-}})]=- 2w
\, Q({\hat a_-}), \label{ezf}
\ee
and the antisymmetric part of the coproduct (\ref{eh}) is reduced to
\bea
&& \hat \delta(Q(\hat a_+))=2 \,Q(\hat \te)\wedge Q(\hat a_+),\nonumber\\
&& \hat \delta(Q(\hat a_-))=-2 \,Q(\hat \te)\wedge Q(\hat a_-),\label{ezg}\\
&& \hat \delta(Q(\hat \te))=0.\nonumber
\eea
It is easy to check that $(sb(2),\hat\delta)$ is the dual Lie bialgebra of
$(iso(1,1),\delta)$. Moreover, $(sb(2),\hat\delta)$ is a coboundary bialgebra
generated by the classical $r$ matrix $\hat r=\frac 1 w Q(\hat a_+)\wedge
Q(\hat \te)$ that verifies the CYBE.

This result can be related to the dual picture developed in the previous
section. Provided that $A_+$, $A_-$ and $H$ are considered as classical
(commutative) coordinates, the matrix $T^Q$ (\ref{ez}) can be viewed as a
$4\times 4$ realization of the group $SB(2)$. Then, we obtain:

\noindent {\bf Proposition 4.4.}
\begin{em}
The fundamental Poisson brackets
\be
\pois{H}{A_+}=\frac{e^{2w{A_+}}-1}w, \quad \pois{H}{A_-}=-
2\, {A_-}e^{2w{A_+}} , \quad  \pois{A_+}{A_-}=0,
\label{ezh}
\ee
endow $Fun(SB(2))$ with a Poisson--Lie structure.
\end{em}

The coboundary structure of $(sb(2),\hat\delta)$ is essential in order to
construct the bracket (\ref{ezh}) by using the Sklyanin procedure (\ref{ek}).
The left and right invariant vector fields
\be
X_{Q(\hat \te)}^L=\partial_H,\qquad
X_{Q(\hat a_+)}^L= 2\,w\,H\,\partial_{H} +\partial_{A_+},\qquad
X_{Q(\hat a_-)}^L=e^{-2wA_+}\partial_{A_-},\label{ezi}
\ee
\be
X_{Q(\hat \te)}^R=e^{-2wA_+}\partial_{H},\qquad
X_{Q(\hat a_+)}^R= - 2\,w\,A_-\,\partial_{A_-} +\partial_{A_+},\qquad
X_{Q(\hat a_-)}^R=\partial_{A_-},\label{ezj}
\ee
are obtained from $T^Q$ (\ref{ez}). The  coproduct  (\ref{df})
is now the classical $SB(2)$ group law.  Moreover, the triangular nature of
$\hat r$ ensures the existence of a $\ast_\h$--product quantizing
(\ref{ezh}). A similar construction has been fully developed in \cite{PL} for
the standard deformation of the Heisenberg algebra.


\sect{A non--standard quantum Heisenberg group}

It is well known that the points of the two dimensional Euclidean plane can
be identified with the {\it complex numbers}, that is, each point with
Cartesian coordinates $(x_1,x_2)$ corresponds to the complex number
$Z=x_1+ix_2$. A similar description can be done for the (1+1) Galilei
and Poincar\'e planes by considering {\it dual numbers} and {\it double
numbers}, respectively. Thus, points within the planes of these three
geometries can be formally described in terms of numbers $Z=x_1+jx_2$,  where
the ``imaginary" unit $j$ is outside the real numbers, and can be equal either
to the complex unit $i$ $(i^2=-1)$, to the dual unit $\epsilon$
$(\epsilon^2=0)$ or to the double unit $e$ $(e^2=1)$ \cite{Yaglom,Grom2}. In
this section we show how it is possible to extend formally the non-standard
quantum Poincar\'e group studied in the previous Section to
non--standard quantum analogues of the (1+1) Euclidean and
non--extended Galilei (Heisenberg) groups. This approach reveals also some
traits of the non--standard Poincar\'e group which are somewhat hidden in
the previous treatment.

If we perform the following formal change of basis
(depending on the unit $j$) in (\ref{bc}):
\be
P_1=\frac 12 (P_- + P_+),\quad P_2=\frac 12 j(P_- -
P_+),\quad J_{12}=-\frac 12 j K ,\label{fc}
\ee
the commutation relations (\ref{bc}) turn into:
\be
[J_{12},P_1]=P_2,\qquad [J_{12},P_2]=j^2 P_1,\qquad [P_1,P_2]=0.\label{fd}
\ee
The change (\ref{fc}) gives rise to a transformation of
the group coordinates. By using the classical exponential factorization we
have that
\bea
&& e^{{ a_-}  P_- }e^{{ a_+}  P_+ }e^{{ \te}   K }=
e^{{ a_-} (P_1+\frac 1j P_2)}e^{{ a_+} (P_1-\frac 1j P_2)} e^{{
\te} (-\frac 2j J_{12})},\nonumber\\
&& \qquad e^{({ a_-}+{ a_+})P_1}e^{\frac 1j({ a_-} - {
a_+})P_2}e^{-\frac 2j{ \te} J_{12}} =e^{a_1P_1}e^{a_2P_2}e^{\de J_{12}}.
\label{ff}
\eea
Hence, we define a quantum analogue of the expressions (\ref{ff})
\be
 \hat a_1={\hat a_-} + {\hat a_+}, \qquad
\hat a_2=\frac 1j ({\hat a_-}  - {\hat a_+}),
\qquad \hat \de=-\frac 2j{\hat \te},
\label{fg}
\ee
 completed with a suitable defined transformation of the
deformation parameter $w$ in the form
\be
\vvv=-\frac 2j w.
\label{fh}
\ee

Strictly speaking, expressions   (\ref{fg}) and (\ref{fh}) are only  defined
for  $j\ne \epsilon$, because  the dual unit has no inverse. However, the
  Hopf algebras obtained by means of these substitutions are always
well--defined:

\noindent {\bf Proposition 5.1.}
\begin{em}
 The  Hopf algebra
$Fun_\vvv(ISO(1,1;j))$  has multiplication, coproduct, counit and antipode
given by
\bea
&& [\hat \de,\hat a_1]=\vvv\left((\cosh j\hat \de - 1) -
j\frac {\sinh j\hat\de}{j}\right),\cr
&& [\hat \de,\hat a_2]=   \vvv \left(
\frac{ \sinh j\hat \de}{j}  -
j\frac{(\cosh j\hat \de -1)}{j^2}\right),\label{fi}\\
&& [\hat a_1,\hat a_2]=
\vvv(\hat a_1+j\hat  a_2),\nonumber
\eea
\bea
&&\Delta(\hat \de)=\hat \de\otimes 1 + 1 \otimes\hat  \de,\cr
&&\Delta(\hat a_1)=\hat a_1\otimes 1 + \cosh j\hat \de\otimes \hat  a_1 + j
\sinh j\hat \de\otimes \hat a_2 , \label{fj}\\
&&\Delta(\hat a_2)=\hat a_2\otimes 1 + \cosh j\hat \de\otimes\hat   a_2 +
\frac { \sinh j\hat \de}j\otimes\hat  a_1 ,\nonumber
\eea
\be
\epsilon(X) =0,\qquad
X\in\{ \hat   a_1,\hat    a_2, \hat \de\};
\label{fk}
\ee
\bea
&&\gamma( \hat \de)=-\hat  \de,\cr
&&\gamma(\hat  a_1)=-\cosh (j\hat \de)\hat  a_1 + j \sinh (j\hat \de) \hat
a_2, \label{fl}\\
&&\gamma(\hat  a_2)=-\cosh (j\hat \de)\hat  a_2 + \frac  {\sinh (j\hat
\de)\hat  a_1}j.\nonumber \eea \end{em}

By taking into account the power series expressions for the  functions
  $\sinh (j\hat y)/j$,
$\cosh (j\hat y)-1$ and $(\cosh
(j\hat y)-1)/j^2$, it is easy to check that all of them  are always real
functions for real values of the argument $y\in \R$,  no matter of  wether $j$
is the complex, dual or double unit. Therefore, note that each commutator
(\ref{fi}) is written as a pair formed by a ``real" term and a ``pure
imaginary" term.

The above statement enables
us to define the quantum Euclidean, Poincar\'e and Galilean planes as follows:

\noindent {\bf Proposition 5.2.}
\begin{em}
Consider the coordinates of the
quantum space $(\hat x_1 ,\hat x_2)$ verifying
\be
[\hat x_1,\hat x_2]= \vvv (\hat x_1 +j \hat x_2),\label{fm}
\ee
then,  the co--action
\be
\left(\begin{array}{c}
1\\ \hat x_1'\\ \hat x_2' \end{array}\right)
=\left(\begin{array}{ccc}
 1 & 0 & 0 \\ \hat a_1 & \cosh j\hat \de & j\sinh j\hat \de\\
\hat a_2 &\frac 1j\sinh j\hat \de & \cosh j\hat  \de  \end{array}\right)
\dot\otimes \left(\begin{array}{c}
1\\ \hat x_1\\ \hat x_2 \end{array}\right),
\label{fn}
\ee
preserves the commutation rules (\ref{fm}) for the $(\hat x'_1 ,\hat x'_2)$
coordinates.
\end{em}

In the Euclidean case $(j=i)$, the resultant Hopf algebra can be
compared to the standard one obtained in \cite{PL} leading to similar
comments to the ones done in $\S$4.1. For the Heisenberg case, the fact that
$j^2=\epsilon^2=0$ implies that commutation rules (\ref{fi}) can be
written as quadratic relations:
\bea
&& [\de,a_1]= - \vvv\, \epsilon\,
\de,\cr && [\de,a_2]=   \vvv \, (  \de - \epsilon \,\de^2/2),\label{fo}\\
&& [a_1,a_2]= \vvv\,(a_1+\epsilon \, a_2). \nonumber
\eea
If we rename generators in the form
$\hat\de=\delta,\hat a_1=\alpha$ and $\hat a_2=\beta$, this algebra is
included within the classification of quantum Heisenberg groups given in
\cite{HLR} (see also \cite{Kup}) as a (Type I) two--parametric deformation
with $z=\frac{1}{2} \vvv \epsilon$ and $p=q=-\vvv$.

It is interesting to note that the standard quantum Heisenberg group
\cite{PL} is obtained by supressing in (\ref{fo}) all ``imaginary" terms that
include the dual unit $\epsilon$. The same is true for the Poincar\'e (resp.
Euclidean) cases if we annihilate the   imaginary
parts of the expressions. In this sense, it is worth remarking that an
approach with dual units is quite different from a strict  contraction
procedure. Thus, if we perform the   contraction \cite{PL}
\be
\hat a'_1=\hat a_1,\qquad \hat a'_2=\lambda^{-1}\hat a_2,\qquad
\hat \de' =\lambda^{-1}\hat \de ,\qquad \vvv'=\lambda^{-1}\vvv,\qquad
(\lambda\to 0),\label{fp}
\ee
in either the non--standard
quantum Poincar\'e or Euclidean groups we get exactly the standard
quantum Heisenberg group.

The geometrical meaning of the relations (\ref{fi})
and (\ref{fm}) should deserve further study; in this non--standard quantum
space the commutator of the coordinates $(\hat x_1,\hat x_2)$ of a  point is,
in some sense, the point itself.


\sect{Concluding remarks}

Non--standard quantum deformations have received scant attention compared to
the standard ones. However, they present interesting features: at a purely
mathematical level, the existence of a $\ast_\h$--product that quantizes
the Poisson--Lie group is always guaranteed for them. From a physical point
of view, they are naturally adapted to the null--plane basis of the
Poincar\'e algebra (see \cite{Bey,NP} for the construction of such null--plane
deformations of the (2+1) and (3+1) dimensional cases).
This paper can be seen as a first step in the study of these
null--plane quantum groups. At this respect, the symmetric appearance  of both
spatial and time quantum translations, together with the coboundary  nature of
the underlying Lie bialgebras are worth emphasizing.

It is also interesting to
note that the essential features of the $T$--matrix approach can be obtained
without computing all the components of the dual tensors. We have  also shown
that the quantum group $Fun_w(ISO(1,1))$ is, in fact, an infinite dimensional
Hopf--Lie algebra. This structure is rather
different to the universal enveloping algebras encountered when computing
$T$--matrices in the literature.

Since the Sklyanin bracket defines a quadratic
algebra in terms of the group entries, it seems difficult, in general,  to
deduce from it the algebraic properties that characterize the dynamical
algebras comming from different $r$--matrices. However, the results contained
in the previous pages show that, if we take into account the specialization of
the Sklyanin bracket to the group coordinates, the comparison between two
Sklyanin algebras defined on the same group is much more easy. The $T$--matrix
provides a canonical framework to consider the corresponding
non--commutative problem and, for the cases so far explored, the main
result  is: at the algebra level, quantum groups are just the Weyl
quantization of these ``fundamental Sklyanin algebras", perhaps endowed
with a deformed group law.

An interesting question to be posed now is whether the study of the
non--standard $sl(2,\R)$ deformation from the $T$--matrix point of view would
confirm the previous analysis.  Work on this line is currently in progress.


\bigskip
\bigskip

\noindent
{\large{{\bf Acknowledgements}}}

\bigskip

C.M.P. expresses her gratitude for the warm hospitality during her stays
in Valladolid. This work has been partially supported by a DGICYT project
(PB92--0255) from the Ministerio de Educaci\'on y Ciencia de Espa\~na.

\bigskip


\end{document}